\documentclass[hyper,11pt]{JHEP3}

\usepackage{xy,amssymb,amsmath,diagrams,times,theorem}

\xyoption{poly}
\xyoption{arc}
\xyoption{knot}
\xyoption{curve}
\xyoption{arrow}
\xyoption{all}

\newcommand{\C}{\mathbb{C}}
\newcommand{\Q}{\mathbb{Q}}

\newcommand{\N}{\mathbb{N}}
\newcommand{\Z}{\mathbb{Z}}
\newcommand{\Ascr}{{\mathcal A}}
\newcommand{\Oscr}{{\mathcal O}}

\newcommand{\V}{\mathbb{V}}
\newcommand{\PP}{\mathbb{P}}
\newcommand{\mumu}{{\pmb \mu}}
\newcommand{\een}{{\xy 0;/r.08pc/:
(3,-3); (5,-3) **@{-}; (5,5) **@{-}; (0,2) **@{-};
(3,-3); (3,4) **@{-}
\endxy}}

\newcommand{\wis}[1]{{\text{\em \usefont{OT1}{cmtt}{m}{n} #1}}}

\newtheorem{theorem}{Theorem}
\newtheorem{lemma}{Lemma}

{\theorembodyfont{\rmfamily} 
\newtheorem{example}{Example}

\newtheorem{definition}{Definition} }

\preprint{UIA preprint 2002-08}

\title{Noncommutative Smooth Models}

\author{Lieven Le Bruyn \\
Departement Wiskunde en Informatica, Universiteit Antwerpen (UIA) \\
B-2610 Antwerp (Belgium) \\
E-mail : \email{lebruyn@uia.ua.ac.be}}

\abstract{We determine the central simple algebras $\Sigma$ over a functionfield $K$ of trancendence degree two which admit a model of smooth Cayley-Hamilton algebras.  This happens
if and only if there is a smooth model $S$ of $K$ such that the ramification divisor of a maximal
$\Oscr_S$-order in $\Sigma$ is a disjoint union of smooth curves. Further, we prove that the
Brauer-Severi fibration of smooth models which are maximal orders is a flat morphism.}

\begin{document}

One can define noncommutative smooth orders either by homological or by
geometrical properties. The first approach has the advantage that there are enough
regular algebras around to have smooth models at least in central simple algebras over 
surfaces, see \cite{ChanIngalls}. The second approach, based on smooth Cayley-Hamilton orders,
is already tedious in dimension two but has the advantage of having an \'etale local description
which makes it possible (at least in principle) to generalize the results in this paper to arbitrary dimensions.

\noindent
{\bf Acknowledgement : }
This paper is a slightly expanded version of a talk given at the 60th birthday conference for
Claudio Procesi, june 2001 in Roma. I like to thank the organizers for the invitation.

\section{Smooth Cayley-Hamilton models}

Let $K$ be a functionfield of trancendence degree $d$ over an algebraically closed field of
characteristic zero which we denote by $\C$. Let $\Sigma$ be a central simple $K$-algebra of
dimension $n^2$ and let $tr$ be the reduced trace of $\Sigma$. Let $C$ be an affine normal domain
with function field $K$ and let $A$ be a $C$-order in $\Sigma$, that is, $A$ is an affine algebra
with center $C$ such that $A.K = \Sigma$. Because $C$ is integrally closed it follows that 
the reduced trace $tr$ on $\Sigma$ makes $A$ into a Cayley-Hamilton algebra of degree $n$, see
\cite{ProcesiCH} or \cite{LBlocalstructure}, such that
$tr(A) = C$. 

Let $\wis{trep}_n~A$ be the affine scheme of all $n$-dimensional {\em trace preserving} representations
of $A$, that is, all $\C$-algebra morphisms $\phi~:~A \rTo M_n(\C)$ compatible with the trace maps
\[
\begin{diagram}
A & \rTo^{\phi} & M_n(\C) \\
\dTo^{tr} & & \dTo^{Tr} \\
A & \rTo^{\phi} & M_n(\C)
\end{diagram}
\]
where $Tr$ is the usual trace on $M_n(\C)$. Clearly, the scheme $\wis{trep}_n~A$ comes equipped with
a natural $GL_n$ (actually $PGL_n$) action and it follows from \cite{Artin} and \cite{ProcesiCH}
that the geometric points of the quotient scheme
\[
\wis{triss}_n~A = \wis{trep}_n~A // GL_n \]
are in one-to-one correspondence with the isomorphism classes of $n$-dimensional semisimple
trace preserving representations of $A$. Moreover, we can recover the algebras $C$ and $A$ from the
$GL_n$-action
\[
C = \C[\wis{triss}_n~A] \qquad A = M_n(\C[\wis{trep}_n~A])^{GL_n} \]
as respectively the algebra of polynomial $GL_n$-invariants on $\wis{trep}_n~A$ and the algebra of 
$GL_n$-equivariant maps from $\wis{trep}_n~A$ to $M_n(\C)$, see \cite{ProcesiCH}.

\begin{definition} With notation as before, $A$ is said to be a {\em smooth} Cayley-Hamilton order if 
and only if $\wis{trep}_n~A$ is a smooth $GL_n$-scheme (in particular reduced). 
\end{definition}

In \cite{ProcesiCH} it
is shown that this implies that $A$ satisfies the extension of Grothendieck's characterization of commutative regular 
algebras (by the universal lifting property modulo nilpotent ideals) to the category of all algebras
satisfying the $n$-th formal Cayley-Hamilton equation. 

\begin{definition} With notations as before, $\Sigma$ is said to possess a {\em noncommutative smooth model}
if there is a projective normal variety $X$ with functionfield $\C(X) = K$ and a sheaf $\Ascr$ of $\Oscr_X$-orders
in $\Sigma$ such that for an affine cover $\{ U_i \}$ of $X$ the sections
\[
A_i = \Gamma(U_i,\Ascr) \]
are smooth Cayley-Hamilton orders. That is, there is a smooth $GL_n$-variety $\wis{trep}_n~\Ascr$
with algebraic quotient $\wis{trep}_n~\Ascr // GL_n = X$.
\end{definition}

If $X$ is a smooth projective model of $K$ and if $\Sigma$ contains a sheaf of Azumaya $\Oscr_X$-algebras
$\Ascr$, then $\Sigma$ has a noncommutative smooth model as $\wis{trep}_n~\Ascr$ is a principal
$PGL_n$-bundle over $X$. Hence, the problem reduces to determining which ramified classes of the Brauer
group contain a noncommutative smooth model.

A natural strategy to prove that every central simple $K$-algebra $\Sigma$ possesses a noncommutative smooth
model would be the following. Let $X$ be a smooth projective model of $K$ and let $\Ascr$ be a sheaf of
maximal $\Oscr_X$-orders in $\Sigma$ and construct the scheme $\wis{trep}_n~\Ascr$. Usually, this scheme
will have singularities but we can construct a $GL_n$-equivariant desingularization of it
\[
\begin{diagram}
\widetilde{\wis{trep}_n~\Ascr} & \rOnto & \wis{trep}_n~\Ascr \\
\dOnto^{\tilde{\pi}} & & \dOnto \\
\tilde{X} & \rOnto & X
\end{diagram}
\]
If the $GL_n$-quotient variety $\tilde{X}$ of this desingularization has an affine open cover
$\{ V_i \}$ such that there are $\Gamma(V_i,\Oscr_{\tilde{X}})$-orders $A_i$ such that
\[
\tilde{\pi}^{-1}(V_i) \simeq \wis{trep}_n~A_i \]
as $GL_n$-varieties, then $\Sigma$ would have a noncommutative smooth model. Unfortunately, this approach
is far too optimistic.

\begin{example}
Let $A$ be the {\em quantum plane} of order two,
\[
A = \frac{\C \langle x,y \rangle}{(xy+yx)} \]
One verifies that $u=x^2$ and $v=y^2$ are
central elements of $A$ and that $A$ is a free module of rank $4$ over $\C[u,v]$. In fact, $A$ is a
$\C[u,v]$-order in the quaternion division algebra 
\[
\Sigma = \begin{pmatrix} u & & v \\
& \C(u,v) & \end{pmatrix}
\]
The induced reduced trace $tr$ is the linear map on $A$ such that 
\[
\begin{cases}
tr(x^iy^j) = 0 &\quad \text{if either $i$ or $j$ are odd, and} \\
tr(x^iy^j) = 2x^iy^j &\quad \text{if $i$ and $j$ are even.}
\end{cases}
\]
In particular, a trace preserving $2$-dimensional representation is determined by a couple of
$2 \times 2$ matrices
\[
(~\begin{bmatrix} x_1 & x_2 \\ x_3 & -x_1 \end{bmatrix}~,~\begin{bmatrix}
x_4 & x_5 \\ x_6 & -x_4 \end{bmatrix}~) \quad \text{with}
\quad~tr(\begin{bmatrix} x_1 & x_2 \\ x_3 & -x_1 \end{bmatrix}~.~\begin{bmatrix}
x_4 & x_5 \\ x_6 & -x_4 \end{bmatrix}) = 0 \]
That is,$\wis{trep}_2~A$ is the hypersurface in $\C^6$ determined by the equation
\[
\wis{trep}_2 A = \mathbb{V}(2x_1x_4+x_2x_6+x_3x_5) \rInto \C^6 \]
and is therefore irreducible of dimension $5$ with an isolated singularity at $p=(0,\hdots,0)$.

Consider the blow-up of $\C^6$ at $p$ which is the closed subvariety of
$\C^6 \times \PP^5$ defined by
\[
\tilde{\C}^6 = \V(x_iX_j-x_jX_i) \]
with the $X_i$ the projective coordinates of $\PP^5$. The strict transform of $\wis{trep}_2~A$
is the subvariety
\[
\widetilde{\wis{trep}_2~A} = \V(x_iX_j-x_jX_i,2X_1X_4+X_2X_6+X_3X_5) \rInto \C^6 \times \PP^5 \]
which is a smooth variety. Moreover, there is a natural $GL_2$-action on it induced by
simultaneous conjugation on the fourtuple of $2 \times 2$ matrices
\[
\begin{bmatrix} x_1 & x_2 \\ x_3 & -x_1 \end{bmatrix} \quad
\begin{bmatrix} x_4 & x_5 \\ x_6 & -x_4 \end{bmatrix} \quad
\begin{bmatrix} X_1 & X_2 \\ X_3 & -X_1 \end{bmatrix} \quad
\begin{bmatrix} X_4 & X_5 \\ X_6 & -X_4 \end{bmatrix} \]
As the projection $\widetilde{\wis{trep}_2~A} \rOnto \wis{trep}_2 A$ is a $GL_2$-isomorphism outside
the exceptional fiber, we only need to investigate the semi-stable points over $p$. Take the
particular point $x$
\[
\begin{bmatrix} 0 & 0 \\ 0 & 0 \end{bmatrix} \quad
\begin{bmatrix} 0 & 0 \\ 0 & 0 \end{bmatrix} \quad
\begin{bmatrix} i & 0 \\ 0 & -i \end{bmatrix} \quad
\begin{bmatrix} 0 & a \\ -a & 0 \end{bmatrix} \]
which is semi-stable and has as stabilizer
\[
Stab(x) = \mumu_2 = \langle \begin{bmatrix} 0 & 1 \\ 1 & 0 \end{bmatrix} \rangle \rInto PGL_2 \]
There is no affine $GL_2$-stable open of $\widetilde{\wis{trep}_2~A}$ containing $x$ 
of the form $\wis{trep}_2~B$ for some order $B$ as this would contradict the fact
that the stabilizer subgroup of any finite dimensional representation is connected. 
\end{example}

Therefore, we need a subtler strategy to investigate the obstruction.

\section{The strategy}

The class $[ \Sigma ] \in Br_n(K) = H^2_{et}(K,\mumu_n)$ in the ($n$-torsion part of the)
Brauer group of $K$ can be found using the {\em coniveau spectral sequence}
\[
E_1^{p,q} = \oplus_{x \in X^{(p)}} H^{q-p}_{et}(\C(x),\mumu_n^{\otimes l-p}) \Rightarrow H^{p+q}_{et}(X,\mumu_n^{\otimes l})
\]
where $X$ is a smooth projective model for $K$ and $X^{(p)}$ is the set of irreducible subvarieties
$x$ of $X$ of codimension $p$ with function field $\C(x)$, see \cite{Grothendieck}. In low dimensions, the resulting
sequence for the Brauer group can be expressed in terms of ramification data of maximal orders in the central simple
algebra, see for example the Artin-Mumford exact sequence for the Brauer group of a smooth surface
\cite{AM}. In case $X$ has singularities one can extend this using the Bloch-Ogus coniveau spectral
sequence, \cite{BlochOgus}.

To apply this result to our problem we need to have some control on the central singularities of
smooth Cayley-Hamilton orders and on the \'etale local structure of the smooth order and thereby on its
ramification locus. Both problems can be solved (at least in small dimensions) and we refer for more
details to the papers \cite{LBlocalstructure} and \cite{SOS}. Here, we merely state the results.

Let $\Ascr$ be a sheaf of smooth Cayley-Hamilton $\Oscr_X$-orders in $\Sigma$ where $X$ is a normal
projective variety with affine open cover $\{ U_a \}$ and corresponding sections
$A_a = \Gamma(U_a,\Ascr)$ and $C_a = \Gamma(U_a,\Oscr_X)$. A point $p \in U_a$ determines an isomorphism class
of a semisimple $n$-dimensional representation of $A_a$ say
\[
V_p = S_1^{\oplus e_1} \oplus \hdots \oplus S_k^{\oplus e_k} \]
where $S_i$ is a simple representation of dimension $d_i$ and occurring in $V_p$ with multiplicity $e_i$, that is
$n = \sum_i d_ie_i$. The space of self-extensions
\begin{eqnarray*}
N_p = Ext^1_{A_a}(V_p,V_p) =& \oplus_{i,j} Ext^1_{A_a}(S_i,S_j)^{\oplus e_ie_j} \\
=& \wis{rep}_{\alpha}~Q
\end{eqnarray*}
is the representation space of $\alpha$-dimensional representations with $\alpha = (e_1,\hdots,e_k)$
of the quiver $Q$ on $k$ vertices $\{ v_1,\hdots,v_k \}$ such that the number of directed arrows
between $v_i$ and $v_j$ is equal to $dim_{\C}~Ext^1_{A_a}(S_i,S_j) - \delta_{ij}$.

The normal space $N_p$ to the orbit of $V_p$ in $\wis{trep}_n~\Ascr$ is the space of trace-preserving self-extensions and is
therefore a linear subspace of $\wis{rep}_{\alpha}~Q$ which can be described as the representation space
of $\alpha$-dimensional representations $\wis{rep}_{\alpha}~Q^{\bullet}$ of a {\em marked} subquiver
$Q^{\bullet}$ of $Q$. That is, some loops of $Q$ may acquire a marking and the matrices corresponding to these
loops have to be of trace zero in $\wis{rep}_{\alpha}~Q^{\bullet}$. The stabilizer subgroup of $V_p$
is the group $GL(\alpha) = GL_{e_1} \times \hdots \times GL_{e_k}$ and the action of it on the
normalspace $N_p$ coincides with the basechange action on $\wis{rep}_{\alpha}~Q^{\bullet}$. The relevance of this
description is that by the Luna slice theorems \cite{Luna} the \'etale local structure of $X$ near $p$
is isomorphic to that of the marked quiver quotient variety
\[
\wis{iss}_{\alpha}~Q^{\bullet} = \wis{rep}_{\alpha}~Q^{\bullet} // GL(\alpha) \]
near the trivial representation $\overline{0}$.
Also, the $GL_n$-structure of $\wis{trep}_n~\Ascr$ near the orbit of $V_p$ is \'etale isomorphic to
that of the associated fiber bundle
\[
GL_n \times^{GL(\alpha)} \wis{rep}_{\alpha}~Q^{\bullet} \]
near the orbit of $\overline{(\een_n,0)}$ where the embedding of $GL(\alpha) \rInto GL_n$ is determined
by the $k$-tuple $\gamma = (d_1,\hdots,d_k)$ of the dimensions of the simple components of $V_p$.

These facts allow us to describe the completion $\hat{\Ascr}_p$ of the stalk of $\Ascr$ in $p$ as well
as $\widehat{\Oscr_{X,p}}$ by walking through the marked quiver setting $(Q^{\bullet},\alpha)$ using the
results of \cite{LBProcesi} and \cite{LBlocalstructure}. The ring of polynomial invariants
$\C[\wis{iss}_{\alpha}~Q^{\bullet}]$ is generated by taking traces of oriented cycles in the quiver
$Q$ and $\widehat{\Oscr_{X,p}}$ is the completion of $\C[\wis{iss}_{\alpha}~Q^{\bullet}]$ at the maximal 
graded ideal. If $M_{ij}$ is the $\C[\wis{iss}_{\alpha}~Q^{\bullet}]$-module of oriented paths from
$v_i$ to $v_j$ in $(Q,\alpha)$, then we have a block-decomposition
\[
\hat{\Ascr_p} \simeq
\begin{bmatrix}
M_{d_1}(M_{11}) & M_{d_1 \times d_2}(M_{12}) & \hdots & M_{d_1 \times d_k}(M_{1k}) \\
M_{d_2 \times d_1}(M_{21}) & M_{d_2}(M_{22}) & \hdots & M_{d_2 \times d_k}(M_{2k}) \\
\vdots & \vdots & \ddots & \vdots \\
M_{d_k \times d_1}(M_{k1}) & M_{d_k \times d_2}(M_{k2}) & \hdots & M_{d_k}(M_{kk})
\end{bmatrix} \otimes \widehat{\Oscr_{X,p}} \]
The fact that $\Ascr$ is an order in a central simple algebra translates to the fact that $\alpha$ must be the dimension vector
of a simple representation of $Q$. If $\chi_Q$ is the Euler form of the quiver $Q$ (that is
$\chi_{ij} = \delta_{ij} - \# \{ \text{arrows from $v_i$ to $v_j$} \}$) then by \cite{LBProcesi} $\beta$ is the dimension vector
of a simple representation of $Q$ iff $Q$ is strongly connected and for all vertex dimension $\delta_i$ we have
\[
\chi_Q(\beta,\delta_i) \leq 0 \qquad \chi_Q(\delta_i,\beta) \leq 0 \]
(unless $Q = \tilde{A_k}$ in which case $\beta = (1,\hdots,1)$). Under this condition the dimension
of the quotient variety $\wis{iss}_{\alpha}~Q^{\bullet}$ (which must be equal to $d$, the dimension of $X$)
is given by
\[
d(\alpha) = 1 - \chi_Q(\alpha,\alpha) - \# \{ \text{marked loops} \} \]
We can also read off the \'etale local structure of the ramification locus of $\Ascr$ in $p$. If
$\beta_1,\hdots,\beta_l$ are dimension vectors of simple representations of $Q$ and if $m_1,\hdots,m_l \in \N$ such that
$\alpha = m_1 \beta_1 + \hdots + m_l \beta_l$ then there is a component of the ramification locus in
$p$ of dimension $d(\beta_1) + \hdots + d(\beta_l)$ and locally isomorphic to
\[
\wis{iss}_{\beta_1}~Q^{\bullet} \times \hdots \times \wis{iss}_{\beta_l}~Q^{\bullet} \]

Working backwards, if we fix the central dimension $d$ and the index $n$ of the central simple algebra,
then we have to determine all marked quiver settings $(Q^{\bullet},\alpha)$ such that
\begin{enumerate}
\item{$d = d(\alpha)$ and $\alpha$ is the dimension vector of a simple representation of $Q^{\bullet}$.}
\item{For $\alpha = (e_1,\hdots,e_k)$ determine all $\gamma = (d_1,\hdots,d_k) \in \N_+^k$ such that
$e_1d_1+ \hdots + e_kd_k = n$.}
\end{enumerate}
Then $\hat{A}_{p}$ must be of the form described above corresponding to such a triple $(Q^{\bullet},\alpha,\gamma)$,
its ramification locus in $p$ only depends on $(Q^{\bullet},\alpha)$ as is the \'etale type of $X$ in
$p$. In \cite{LBlocalstructure} a method is given to classify all relevant triples $(Q^{\bullet},\alpha,\gamma)$ in
low dimensions ($\leq 4$) and in \cite{SOS} reduction steps were given to classify all possible central singularities
in low dimensions ($ \leq 6$). As for the singularities, there are none in dimension $\leq 2$ and
in dimensions $3$ (resp. $4,5,6$) there are precisely $1$ (resp. $3,10,53$) types of possible
singularities for smooth Cayley-Hamilton orders.

These explicit \'etale local data combined with the Bloch-Ogus spectral sequence description of the
$n$ part of the Brauer group should be enough information to determine those central simple algebras
$\Sigma$ over $K$ allowing a noncommutative smooth model in low central dimensions. In the next section we will perform the required
calculations when the central dimension $d=2$.

\section{The case of surfaces}

First, we recall the classification of all admissible triples $(Q^{\bullet},\alpha,\gamma)$ for $d=2$
given in \cite[\S 6]{LBlocalstructure}. All admissible (marked) quiver settings $(Q^{\bullet},\alpha)$
are of the form $A_{klm}$, $k,l \geq 0,m > 0$
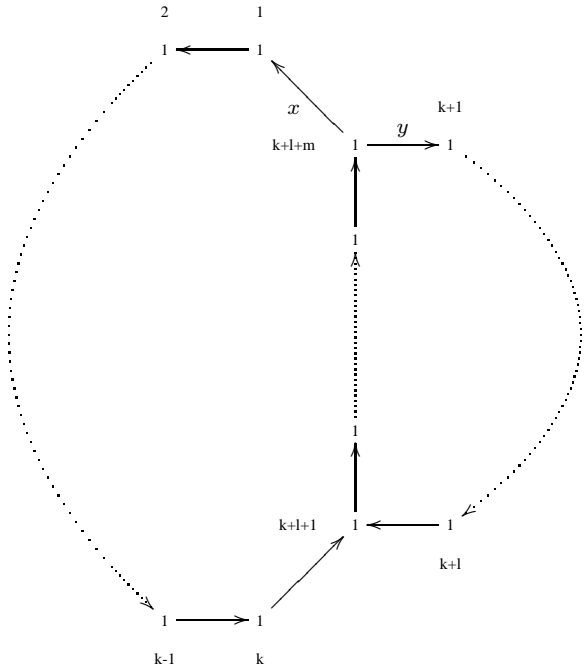
\[ 
\begin{xy}/r.15pc/:
\POS (0,0) *+{\txt{\tiny 1}} ="a" , (20,0) *+{\txt{\tiny 1}}="b",
(40,20) *+{\txt{\tiny 1}}="c",
(40,40) *+{\txt{\tiny 1}}="d",
(40,80) *+{\txt{\tiny 1}}="e",
(40,100) *+{\txt{\tiny 1}}="f",
(20,120) *+{\txt{\tiny 1}}="g",
(0,120) *+{\txt{\tiny 1}}="h",
(60,100) *+{\txt{\tiny 1}}="i",
(60,20) *+{\txt{\tiny 1}}="j"
\POS"a" \ar  "b"
\POS"b" \ar "c"
\POS"c" \ar "d"
\POS"e" \ar  "f"
\POS"f" \ar^{x}  "g"
\POS"g" \ar "h"
\POS"f" \ar^{y}  "i"
\POS"j" \ar "c"
\POS"d" \ar @{.>} "e"
\POS"h" \ar @{.>}@/_12ex/ "a"
\POS"i" \ar @{.>}@/^10ex/ "j"
\POS"g"+(0,8) *+\txt{\tiny 1}
\POS"b"-(0,8) *+\txt{\tiny k}
\POS"a"-(0,8) *+\txt{\tiny k-1}
\POS"h"+(0,8) *+\txt{\tiny 2}
\POS"i"+(0,8) *+\txt{\tiny k+1}
\POS"j"-(0,8) *+\txt{\tiny k+l}
\POS"c"-(12,0) *+\txt{\tiny k+l+1}
\POS"f"-(13,0) *+\txt{\tiny \txt{\tiny k+l+m}}
\end{xy}
\]
where we make the obvious changes whenever $k$ and/or $l$ is zero. If $p=k+l+m$ then the admissible
vectors $\gamma$ are the unordered partitions $\gamma = (d_1,\hdots,d_p)$ of $n$ having exactly $p$
parts. This allows us to determine all possible \'etale local structures of smooth Cayley-Hamilton
$\Oscr_S$-orders $\Ascr$ in $\Sigma$ where $S$ is a normal projective surface with function field $K$.
If $p \in S$ has corresponding data
\[
(A_{klm},(1,\hdots,1),(d_1,\hdots,d_p)) \]
then using \cite[Prop. 6.4]{LBlocalstructure}  with the indicated arrows $x$ and $y$ (one can always reduce 
to this case using the $GL(\alpha)$-action) we have
\[
\widehat{\Oscr_{S,p}} \simeq \C [[x,y]] \]

\par \vskip 4mm

\[
\hat{A}_p \simeq
\setlength{\unitlength}{.50mm}
\begin{picture}(100,50)(0,50)
\put(0,0){\line(0,1){100}}
\put(100,0){\line(0,1){100}}
\put(33,2){\line(0,1){96}}
\put(66,2){\line(0,1){96}}
\put(2,33){\line(1,0){96}}
\put(2,66){\line(1,0){96}}
\put(2,98){\line(1,-1){96}}
\put(16,16){\makebox(0,0){$(x)$}}
\put(16,50){\makebox(0,0){$(x)$}}
\put(50,16){\makebox(0,0){$(y)$}}
\put(50,84){\makebox(0,0){$(y)$}}
\put(84,50){\makebox(0,0){$(1)$}}
\put(84,84){\makebox(0,0){$(1)$}}
\put(24,92){\makebox(0,0){$(1)$}}
\put(58,58){\makebox(0,0){$(1)$}}
\put(92,24){\makebox(0,0){$(1)$}}
\put(8,76){\makebox(0,0){$(x)$}}
\put(42,42){\makebox(0,0){$(y)$}}
\put(76,8){\makebox(0,0){$(x,y)$}}
\put(16,-5){\makebox(0,0){$\underbrace{\quad \quad \quad \quad }_k$}}
\put(50,-5){\makebox(0,0){$\underbrace{\quad \quad \quad \quad }_l$}}
\put(84,-5){\makebox(0,0){$\underbrace{\quad \quad \quad \quad }_m$}}
\end{picture}
\rInto M_n(\C[[x,y]]) \]

\par \vskip 35mm
\noindent
where at place $(i,j)$ for $1 \leq i,j \leq p$ there is a block of dimension $d_i \times d_j$ with entries
in the indicated ideal of $\C[[x,y]]$.
Some immediate consequences can be drawn from these descriptions :
\begin{enumerate}
\item{$p$ is a smooth point of $S$. Hence, the central variety $S$ must be a smooth projective surface.}
\item{$\Ascr$ is \'etale splittable in $p$ as $\hat{\Ascr}_p$ is an order in $M_n(\C[[x,y]])$.}
\item{For $A_{001}$, $\Ascr_p$ is an Azumaya algebra. In all other cases, $p$ is a point of the ramification locus.
There are three possible subcases :
\begin{itemize}
\item{$p$ is an isolated point of the ramification locus in case $A_{00m}$ with $m > 1$.}
\item{$p$ is a smooth point of a one-dimensional branch of the ramification locus through $p$ in
case $A_{k01}$ or $A_{0l1}$.}
\item{The ramification locus has a normal crossing at $p$ in all other cases.}
\end{itemize}}
\end{enumerate}

Therefore we may assume that the central variety $S$ of a smooth Cayley-Hamilton $\Oscr_S$-order in
$\Sigma$ is a smooth projective variety. In this case the coniveau spectral sequence describing the
Brauer group $Br(K)$ is known as the Artin-Mumford exact sequence \cite{AM}. Consider the sequence of
Abelian groups
\[
0 \rTo Br(S) \rTo^i Br(K) \rTo^a \oplus_C H^1(\C(C),\Q/\Z) \rTo^r \oplus_p \mumu^{-1} \rTo^s \mumu^{-1} \rTo 0
\]
where the first sum runs over all irreducible curves $C$ in $S$ and the second over all points
$p \in S$ and
where the maps have the following ringtheoretic interpretation (see \cite{AM} or \cite{Tannenbaum}
for more details). The inclusion $i$ is induced by assigning to an Azumaya $\Oscr_S$-order $\Ascr$ its
central simple ring of central quotients. We have already seen that central simple algebras living in this image
do have noncommutative smooth models.

The map $a$ is determined by taking the ramification divisor of a maximal $\Oscr_S$-order $\Ascr$ in
$\Sigma$. Let $C$ be an irreducible curve in $S$ then it determines a discrete valuation ring $R$ in $K$
with residue field the functionfield $\C(C)$ which has a trivial Brauer group by Tsen's theorem. But then,
\[
\frac{\Ascr \otimes R}{rad~\Ascr \otimes R} \simeq M_s(L) \]
where $L$ is a cyclic field extension of $\C(C)$. The map $a$ sends the class $[ \Sigma ]$ to the class
$[ L ]$ in $H^1(\C(C),\Q/\Z)$ which classifies all cyclic extensions of $\C(C)$. The union of all
curves $C$ for which the extension is non-trivial is said to be the ramification divisor of $\Ascr$
and because $\Ascr$ is a maximal order it coincides with the non-Azumaya locus of $\Ascr$.

The group $\mumu^{-1} = \cup_n~Hom(\mumu_n,\Q/\Z)$ and the map $r$ measures the ramification of the
cyclic extension $L$ of $\C(C)$ in points $c \in C$, or equivalently, of the ramified cyclic covering
of the normalization $\tilde{C}$ of $C$ in all points $\tilde{c}$ of $\tilde{C}$ lying over $c$, see
\cite[p.110]{Tannenbaum}. Finally, $s$ is the sum map.

The Artin-Mumford theorem asserts that the maps $i$ and $a$ form an exact sequence whenever $S$ is a
smooth projective surface and  the full sequence is exact whenever $S$ is in addition simply connected.
We will use a local version of the Artin-Mumford sequence, thus avoiding the issue of simple connectivity.

It will turn out that the requirement of being \'etale splittable at all points is the crucial
obstruction to having a noncommutative smooth model. Therefore, we may restrict attention to 
maximal orders for if $\Ascr$ is a smooth Cayley-Hamilton $\Oscr_S$-order in $\Sigma$ which is
everywhere \'etale splittable, so is every maximal order containing it.

\par \vskip 3mm

Start with a smooth model $S$ of $K$ and take $\Ascr_0$ to be a maximal $\Oscr_S$-order in $\Sigma$
having $D_0$ as ramification divisor which may be highly singular. Using embedded
resolution of singularities (see for example \cite{Hartshorne}) we can construct another smooth
model $\tilde{S}$ of $K$ such that
\[
\begin{diagram}
& & \tilde{S} \\
& & \dOnto_{\pi} \\
D_0 & \rInto & S
\end{diagram}
\]
the inverse image $\pi^{-1}(D_0) = \tilde{D_0} \cup E$ (where $\tilde{D_0}$ is the strict transform
and $E$ is the exceptional fiber) has at worst normal crossings. 

Next, take a maximal $\Oscr_{\tilde{S}}$-order $\Ascr$ in $\Sigma$. Its ramification divisor
$D$ is a subdivisor of $\pi^{-1}(D_0)$ which (as we will see below) may contain components of the
exceptional divisor $E$. Still, $D$ has at worst normal crossings as singularities and we will
determine the \'etale local structure of $\Ascr$ in all points $p \in D$ to verify where it is a
smooth Cayley-Hamilton order. We separate two cases :

\par \vskip 3mm
\noindent
{\bf case 1 : } Let $p$ be a smooth point of $D$. For $U = \wis{spec}~\C[[x,y]]-V(x)$ we have that
$Br_n~U = 0$ and as this is the \'etale local structure of $\tilde{S}$ and $D$ near $p$ it follows that
$\Ascr$ is \'etale splittable in $p$. Moreover, M. Artin determined in \cite{Artinlocal} the \'etale
local structure of maximal orders over surfaces in a smooth point of the ramification divisor. His results
assert that $\hat{\Ascr}_p$ has the block-decomposition
\[
\hat{\Ascr}_p \simeq
\begin{bmatrix}
(1) & \hdots & \hdots & (1) \\
(x) & \ddots &  & \vdots \\
\vdots & \ddots & \ddots & \vdots \\
(x) & \hdots & (x) & (1)
\end{bmatrix}
\]
where $n=a.b$ and all $b^2$ blocks have sizes $a \times a$ and entries in the indicated ideal
of $\C[[x,y]]$. But this \'etale local structure corresponds to the triple
\[
(A_{b-1 01},(1,\hdots,1),(a,\hdots,a)) \]
whence $\Ascr$ is a smooth Cayley-Hamilton order in $p$.

\vskip 3mm
\noindent
{\bf case 2 : } Let $p$ be a normal crossing of the ramification divisor $D$ of $\Ascr$. Let
$U = \wis{spec}~\C[[x,y]] - V(xy)$ then the coniveau spectral sequence (or see \cite{Artinmaxorders})
gives the exact sequence
\[
0 \rTo Br_n~U \rTo \Z_n \oplus \Z_n \rTo^+ \Z_n \rTo 0 \]
and all classes are determined by quantum planes $\C_q[[u,v]]$ where $vu = quv$ where $q$ is an $n$-th root
of one. As this is the \'etale local structure of $\tilde{S}$ and $D$ near $p$, a combination with
the results from \cite{Artinlocal} (or see \cite{ChanIngalls}) asserts that
\[
\hat{A}_p \simeq \begin{bmatrix}
M_c(\C_q [[u,v]] ) & \hdots & \hdots & M_c(\C_q[[u,v]] ) \\
M_c(u \C_q[[u,v]] ) & \ddots &  & \vdots \\
\vdots & \ddots & \ddots & \vdots \\
M_c(u \C_q[[u,v]]) & \hdots & M_c(u \C_q[[u,v]] ) &M_c( \C_q[[u,v]])
\end{bmatrix} \rInto M_{ac}(\C_q[[u,v]])
\]
where $q$ is a primitive $b$-th root of one and $n=a.b.c$. This class corresponds to the local ramification
data
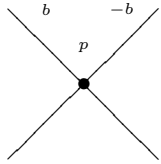
\[
\begin{xy}
(-10,10);(10,-10) **@{-};
(-10,-10);(10,10) **@{-};
(0,0) *+{\bullet};
(0,5) *+{\txt{\tiny{$p$}}};
(-5,10) *+{\txt{\tiny{$b$}}};
(5,10) *+{\txt{\tiny{$-b$}}};
\end{xy}
\]
where the two branches are given classes $b$ and $-b$ in $\Z_n$. As a consequence, $\Ascr$ is not
\'etale splittable in $p$ whenever $b \not= 0$. This problem cannot be resolved by blowing up the point
$p$ as there is only one possibility to satisfy the local sum zero condition for the divisor
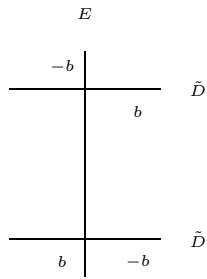
\[
\begin{xy}
(50,-15);(50,15) **@{-};
(40,10);(60,10) **@{-};
(40,-10);(60,-10) **@{-};
(65,10) *+{\txt{\tiny{$\tilde{D}$}}};
(65,-10) *+{\txt{\tiny{$\tilde{D}$}}};
(50,20) *+{\txt{\tiny{$E$}}};
(57,7) *+{\txt{\tiny{$b$}}};
(57,-13) *+{\txt{\tiny{$-b$}}};
(47,13) *+{\txt{\tiny{$-b$}}};
(47,-13) *+{\txt{\tiny{$b$}}};
\end{xy}
\]
which means that the exceptional divisor $E$ will be part of the ramification divisor of a 
maximal order in $\Sigma$ over the blow-up.

\par \vskip 3mm

This concludes the proof of our main theorem :

\begin{theorem} Let $K$ be a functionfield of trancendence degree two over an algebraically closed
field of characteristic zero and let $\Sigma$ be an $n^2$-dimensional central simple $K$-algebra.
Then, $\Sigma$ contains a noncommutative smooth model if and only if there exists a smooth projective
surface $S$ with $\C(S) = K$ such that the ramification divisor of a maximal $\Oscr_S$-order
$\Ascr$ in $\Sigma$ consists of a disjoint union os smooth irreducible curves.
\end{theorem}

If $K = \C(x,y)$ then the $3$-dimensional Sklyanin algebra (see for example \cite{Artinsome}) defines a 
maximal $\Oscr_{\PP^2}$-orders having as ramification divisor a smooth elliptic curve. As a
consequence it is a noncommutative smooth model of the corresponding division algebra. On the contrary,
the triangle $3$-dimensional Auslander regular algebra defines a maximal $\Oscr_{\PP^2}$-order with
ramification divisor three projective lines intersecting transversally but with non-trivial ramification
data. Therefore, it cannot be a noncommutative smooth model for the division algebra of the quantum plane.
In fact, this division algebra cannot contain any noncommutative smooth model.

Another class of division algebras over $\C(x,y)$ possessing a noncommutative smooth model are the
quaternion algebras from \cite{AM} in their construction of unirational non-rational threefolds.
These threefolds are the Brauer-Severi varieties of the corresponding maximal orders over a blow-up.

\section{Brauer-Severi varieties}

Assume we are in the setting of section~1, that is, $C$ is a normal affine domain with field of fractions
$K$ and $A$ is a $C$-order in a central simple $K$-algebra $\Sigma$ of dimension $n^2$. We will define
the {\em Brauer-Severi scheme} $\wis{BSev}~A$ following the account of M. Van den Bergh in \cite{VdBBS}.

There is a natural $GL_n$-action on the product $\wis{trep}_n~A \times \C^n$ defined by
\[
g.(\phi,v) = (g \phi g^{-1}, gv) \]
With $\wis{brauer}~A$ we will denote the open subset of {\em Brauer stable} points of this action
\[
\wis{brauer}~A = \{ (\phi,v)~|~\phi(A) v = \C^n \} \]
which is also the set of points with trivial stabilizer subgroup. Hence, every $GL_n$-orbit in
$\wis{brauer}~A$ is closed and we can form the orbitspace 
\[
\wis{BSev}~A = \wis{brauer}~A / GL_n  \rOnto \wis{triss}_n~A \]
which is known to be a projective space bundle over the quotient variety $\wis{triss}_n~A$(which has
coordinate ring $C$). In general not much can be said about these Brauer-Severi schemes. 

\begin{lemma} If $A$ is a smooth Cayley-Hamilton $C$-order in $\Sigma$, then its Brauer-Severi scheme
$\wis{BSev}~A$ is a smooth variety.
\end{lemma}

\Proof
Because $\wis{brauer}~A \rOnto \wis{BSev}~A$ is a principal $GL_n$-bundle it suffices to show that the
total space $\wis{brauer}~A$ is smooth. This is clear as it is an open subscheme of the smooth variety
$\wis{trep}_n~A \times \C^n$.

\par \vskip 3mm

If $X$ is a normal projective variety with function field $\C(X) = K$ and $\Ascr$ is a $\Oscr_X$-order
in $\Sigma$ we can sheafify the previous construction to obtain a projective space bundle
\[
\wis{BSev}~\Ascr \rOnto^{\pi} X \]
with generic fibers isomorphic to $\PP^{n-1}$ (in the Azumaya points)
and if $\Ascr$ is a noncommutative smooth model, $\wis{BSev}~\Ascr$ is a smooth variety. An important
problem is to determine in which cases $\pi$ is a flat morphism and to determine the structure of all
the fibers $\pi^{-1}(p)$ for $p \in X$.

If $\Ascr$ is a noncommutative smooth model of $\Sigma$, then the \'etale local structure of $\Ascr$
and of $X$ near $p$ is fully determined by a triple $(Q^{\bullet},\alpha,\gamma)$. This triple also
contains enough information to describe the \'etale local structure of the fibration
$\wis{BSev}~A \rOnto X$ near $p$ as the moduli space of a certain quiver setting which we now recall,
more information can be found in \cite{LBSeelinger}.

If $Q^{\bullet}$ is a marked quiver on the vertices $\{ v_1,\hdots,v_p \}$ then we denote with
$\tilde{Q}$ the extended marked quiver obtained from $Q^{\bullet}$ by adding an additional vertex $v_0$
and if $\gamma = (d_1,\hdots,d_p)$ by adding for each $1 \leq i \leq p$ exactly $d_i$ directed arrows
from $v_0$ to $v_i$. If $\alpha = (e_1,\hdots,e_p)$ then we will denote with $\tilde{\alpha}$ the
extended dimension vector $(1,e_1,\hdots,e_p)$. Any character of
$GL(\tilde{\alpha}) = \C^* \times GL(\alpha) \rTo \C^*$ is of the form 
$(g_0,g_1,\hdots,g_p) \mapsto det(g_0)^{t_0}det(g_1)^{t_1}\hdots det(g_p)^{t_p}$ and is therefore
fully determined by the $p+1$-tuple of integers $\theta=(t_0,t_1,\hdots,t_p)$. Consider the character determined
by
\[
\theta = (-n,d_1,\hdots,d_p) \]
then clearly $\theta.\tilde{\alpha} = 0$. In the quiver representation space $\wis{rep}_{\tilde{\alpha}}~\tilde{Q}$
consider the open subset $\wis{rep}^{\theta}$ of $\theta$-semistable representations as defined in
\cite{King}, that is, all representations $V$ such that for every proper subrepresentation $W$ of $V$
we have that $\theta.\wis{dim}(W) \geq 0$. In this particular case we have that all $\theta$-semistable
representations are actually $\theta$-stable, that is, for all proper subrepresentations $W$ of $V$
we have $\theta.\wis{dim}(W) > 0$. The corresponding moduli space of quiver representations
\[
\wis{rep}^{\theta}/GL(\tilde{\alpha}) = \wis{moduli}_{\tilde{\alpha}}^{\theta}~\tilde{Q} \rOnto \wis{iss}_{\tilde{\alpha}}~\tilde{Q} = \wis{iss}_{\alpha}~Q
\]
is a projective space bundle over $\wis{iss}_{\alpha}~Q$ and this bundle is \'etale isomorphic to
the Brauer-Severi fibration in a neighborhood of $p$.

The upshot of this description is that it allows us to calculate the fibers $\pi^{-1}(p)$ of the Brauer-Severi
fibration using the isomorphism (in the Zariski topology)
\[
\pi^{-1}(p) = (\wis{null}_{\tilde{\alpha}}~\tilde{Q} \cap \wis{rep}^{\theta}) / GL(\tilde{\alpha}) \]
where $\wis{null}_{\tilde{\alpha}}~\tilde{Q}$ is the nullcone of quiver representation which can be
analyzed by means of the Hesselink stratification as described in \cite{LBoptimal}.

\par \vskip 3mm

As an example, let us perform the necessary computations in the case of surfaces and a point $p \in S$
where the triple is of the form
\[
(Q^{\bullet} = A_{k01},\alpha = (1,\hdots,1),\gamma = (d_1,\hdots,d_{k+1}) ) \]
and $\gamma$ is a unordered partition of $n$ having precisely $k+1$ parts. Observe that the
\'etale local description of M. Artin of maximal orders over surfaces in smooth points of the ramification
divisor is a special case of such a setting having the additional restriction that $d_1=\hdots=d_{k+1}$.
The extended quiver setting $(\tilde{Q},\tilde{\alpha})$ is of the form
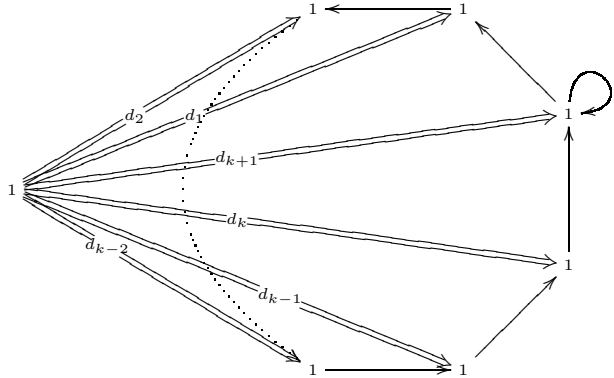
\[
\xy
\POS (0,0) *+{\txt{\tiny $1$}} ="a",
(20,0) *+{\txt{\tiny $1$}} ="b",
(34,14) *+{\txt{\tiny $1$}} ="c",
(34,34) *+{\txt{\tiny $1$}} ="d",
(20,48) *+{\txt{\tiny $1$}} ="e",
(0,48) *+{\txt{\tiny $1$}} ="f",
(-40,24) *+{\txt{\tiny $1$}} ="g"
\POS"a" \ar "b"
\POS"b" \ar "c"
\POS"c" \ar "d"
\POS"d" \ar@(u,r)
\POS"d" \ar "e"
\POS"e" \ar "f"
\POS"f" \ar@/_10ex/@{.>} "a"
\POS"g"  \ar@{=>}|<>(.4)*{\txt{\tiny{$d_{k+1}$}}} "d"
\POS"g"  \ar@{=>}|<>(.4)*{\txt{\tiny{$d_1$}}} "e"
\POS"g"  \ar@{=>}|<>(.4)*{\txt{\tiny{$d_2$}}} "f"
\POS"g"  \ar@{=>}|<>(.3)*{\txt{\tiny{$d_{k-2}$}}} "a"
\POS"g"  \ar@{=>}|<>(.6)*{\txt{\tiny{$d_{k-1}$}}} "b"
\POS"g"  \ar@{=>}|<>(.4)*{\txt{\tiny{$d_k$}}} "c"
\endxy 
\]
and the character corresponds to $\theta = (-n,d_1,\hdots,d_{k+1})$. The weights of the maximal torus
$T = GL(\tilde{\alpha})$ occurring in the representation space $\wis{rep}_{\tilde{\alpha}}~\tilde{Q}$ form the set
\[
\{ \pi_{0i},\pi_{i i+1}~|~1 \leq i \leq k+1 \}
\]
and we know from \cite{LBoptimal} that strata in the nullcone $\wis{null}_{\tilde{\alpha}}~\tilde{Q}$
are determined by saturated subsets of this set of weights. Potential maximal strata correspond to
maximal saturated sets which in this case are exactly the following subsets
\[
S_i = \{ \pi_{i i+1},\pi_{i+1 i+2},\hdots,\pi_{i-2 i-1} \} \cup \{ \pi_{0j}~|~1 \leq j \leq k+1 \}
\]
for every $1 \leq i \leq k+1$. To determine whether the corresponding Hesselink strata of the nullcone
is non-empty we have to determine whether the associated level quiver has $\theta_i$-semistable representations
for a specific character $\theta_i$ (we refer for all these notions to the paper \cite{LBoptimal}).
Allow us to state that the data associated to $S_i$ are : the level quiver data $(Q_i,\alpha_i)$ is
of the form
\[
\begin{xy}
\POS (0,0) *+{\txt{\tiny{$1$}}} ="a0",
(20,0) *+{\txt{\tiny{$1$}}} ="a1",
(40,0) *+{\txt{\tiny{$1$}}} ="a2",
(80,0) *+{\txt{\tiny{$1$}}} ="aE1",
(100,0) *+{\txt{\tiny{$1$}}} ="aE",
(60,0) *+{\txt{$\hdots~\hdots$}} ="dots"
\POS"a0" \ar@{=>}|<>(.4)*{\txt{\tiny{$d_i$}}} "a1"
\POS"a1" \ar "a2"
\POS"a2" \ar "dots"
\POS"dots" \ar "aE1"
\POS"aE1" \ar "aE"
\POS"a0"+(0,5) *+{\txt{\tiny{$-k-1$}}}
\POS"a1"+(0,5) *+{\txt{\tiny{$-k+1$}}}
\POS"a2"+(0,5) *+{\txt{\tiny{$-k+3$}}}
\POS"aE1"+(0,5) *+{\txt{\tiny{$k-1$}}}
\POS"aE"+(0,5) *+{\txt{\tiny{$k+1$}}}
\end{xy} \]
where the superscripts indicate the entries of the character $\theta_i$. Clearly, there are $\theta_i$-semistable
representations in $\wis{rep}_{\alpha_i}~Q_i$ (for example, all arrows equal to $1$) and in fact the
corresponding moduli space is
\[
\wis{moduli}_{\alpha_i}^{\theta_i}~Q_i \simeq \PP^{d_i-1} \]
Hence, the Hesselink strata determined by $S_i$ is non-empty and using the description of
\cite{LBoptimal} the full stratum consists of all representations $\wis{rep}^{ss}_i$ of the extended quiver $\tilde{Q}_i$
which is of the form
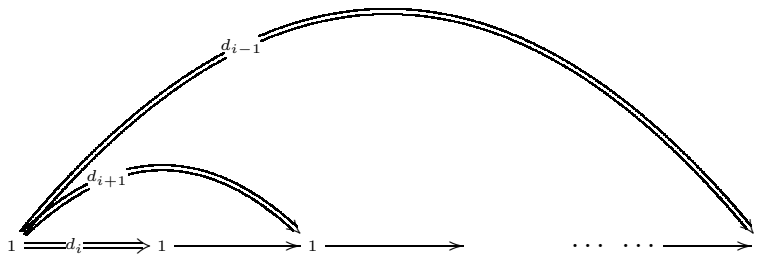
\[
\begin{xy}
\POS (0,0) *+{\txt{\tiny{$1$}}} ="a0",
(20,0) *+{\txt{\tiny{$1$}}} ="a1",
(40,0) *+{\txt{\tiny{$1$}}} ="a2",
(100,0) *+{\txt{\tiny{$1$}}} ="aE",
(60,0) ="dots1",
(80,0) *+{\txt{$\hdots~\hdots$}} ="dots"
\POS"a0" \ar@{=>}|<>(.4)*{\txt{\tiny{$d_i$}}} "a1"
\POS"a1" \ar "a2"
\POS"a2" \ar "dots1"
\POS"dots" \ar "aE"
\POS"a0" \ar@{=>}@/^6ex/|<>(.3)*{\txt{\tiny{$d_{i+1}$}}} "a2"
\POS"a0" \ar@{=>}@/^18ex/|<>(.3)*{\txt{\tiny{$d_{i-1}$}}} "aE"
\end{xy} \]
such that the part of the representation in the level quiver $Q_i$ is $\theta_i$-semistable. That is,
there is an irreducible component of $\wis{null}_{\tilde{\alpha}}~\tilde{Q}$ corresponding to each
$1 \leq i \leq k+1$ each being the closure $\wis{rep}^{ss}_i$ which is of dimension $n+k$. As every point
in the (non-empty) open subset
\[
\wis{rep}^{ss}_i \cap \wis{rep}^{\theta} \]
is $\theta$-stable and hence has stabilizer subgroup $\C^*$ in $GL(\tilde{\alpha})$ which is a
$k+2$-dimensional torus we have that the dimension of the irreducible component in the Brauer-Severi
fiber $\pi^{-1}(p)$ determined by $S_i$ has dimension $n+k-(k+1) = n-1$ for each of the $k+1$ choices
of $i$. Clearly, the fibers of the Brauer-Severi fibration $\wis{BSev}~A \rOnto^{\pi} S$ over an
Azumaya point $p \in S$ are isomorphic to $\PP^{n-1}$ and hence are also of dimension $n-1$. This
completes the proof of

\begin{theorem} Let $\Ascr$ be a noncommutative smooth model of $\Sigma$ over a smooth surface $S$
such that the ramification locus of $\Ascr$ consists (at worst) of a disjoint union of smooth irreducible
curves in $S$. Then, the Brauer-Severi fibration
\[
\wis{BSev}~\Ascr \rOnto^{\pi} S \]
is a flat morphism and the fiber over a ramified point with corresponding marked quiver $A_{k01}$
has exactly $k+1$ irreducible components (each of dimension $n-1$).
\end{theorem}

Recall that the condition is automatic for $\Ascr$ a smooth maximal $\Oscr_S$-order hence this result extends
the calculation of Artin and Mumford in \cite{AM} for the Brauer-Severi scheme of smooth maximal orders 
in specific quaternion algebras.

\providecommand{\href}[2]{#2}
\begingroup\raggedright\endgroup


\begin{thebibliography}{10}



\bibitem{Artin}
Michael Artin, {\it On Azumaya algebras and finite dimensional representations of rings} J. Alg.
{\bf 11} (1969) 523-563

\bibitem{Artinlocal}
Michael Artin, {\it Local structure of maximal orders on surfaces} in "Brauer groups in ring theory
and algebraic geometry" (Antwerp 1981) Lect. Notes Math. 917 (1982) Springer-Verlag

\bibitem{Artinmaxorders}
Michael Artin, {\it Maximal orders of global dimension and Krull dimension two}, Invent. Math. {\bf 84} (1986) 195-222

\bibitem{Artinsome}
Michael Artin, {\it Some problems on three-dimensional graded domains} in "Representation theory
and algebraic geometry" (Waltham 1995) LMS Lect. Notes {\bf 238} (1997) 1-19

\bibitem{AM}
Michael Artin and David Mumford, {\it Some elementary examples of unirational varieties
which are not rational}, Proc. LMS {\bf 25} (1972) 75-95

\bibitem{BlochOgus}
S. Bloch and A. Ogus, {\it Gersten's conjecture and the homology of schemes}
Ann. Sci. ENS {\bf 7} (1974) 181-202

\bibitem{SOS}
Raf Bocklandt, Lieven Le Bruyn and Geert Van de Weyer, {\it Smooth order singularities}
[\href{http://xxx.lanl.gov/abs/math.RA/0207250}{{\tt arXiv:math.RA/0207251}}] (to appear)

\bibitem{ChanIngalls}
Daniel Chan and Colin Ingalls, {\it Birational classification of orders over surfaces} preprint
(2002) to appear

\bibitem{Grothendieck}
Alexander Grothendieck, {\it Le groupe de Brauer I,II,III} in "Dix expos\'es sur la cohomologie des sch\'emas",
North-Holland (1968) 46-188

\bibitem{Hartshorne}
Robin Hartshorne, {\it Algebraic Geometry} GTM {\bf 52} (1977) Springer-Verlag

\bibitem{King}
Alastair King, {\it Moduli of representations of finite dimensional algebras}, Quat. J. Math. Oxford
{\bf 45} (1994) 515-530


\bibitem{LBlocalstructure}
Lieven Le Bruyn, {\it Local structure of Schelter-Procesi smooth orders}, Trans. AMS {\bf 352}
(2000) 4815-4841

\bibitem{LBoptimal}
Lieven Le Bruyn, {\it Optimal filtrations on representations of finite dimensional algebras},
Trans. AMS {\bf 353} (2001) 411-426

\bibitem{LBProcesi}
Lieven Le Bruyn and Claudio Procesi, {\it Semi-simple representations of quivers},
Trans AMS {\bf 317} (1990) 585-598

\bibitem{LBSeelinger}
Lieven Le Bruyn and George Seelinger, {\it Fibers of generic Brauer-Severi schemes},
J. Alg. {\bf } 214 (1999) 222-234.
\bibitem{Luna}
Domingo Luna, {\it Slices \'etales} Bull. Soc. Math. France, M\'emoire {\bf 33} (1973) 81-105

\bibitem{ProcesiCH}
Claudio Procesi, {\it A formal inverse to the Cayley-Hamilton theorem}
J.Alg. {\bf 107} (1987) 63-74



\bibitem{Tannenbaum}
Allen Tannenbaum, {\it The Brauer group and unirationality : an example of Artin-Mumford}
in "Groupe de Brauer" (Les Plans-sur-Bex 1980) Lect. Notes Math. 844 (1981) Springer-Verlag

\bibitem{VdBBS}
Michel Van den Bergh, {\it The Brauer-Severi scheme of the trace ring of generic matrices}
in "Perspectives in Ring Theory" (Antwerp 1986) Kluwer (1988)



\end{thebibliography}
\end{document}